\documentclass[a4paper]{article}
\usepackage[a4paper]{geometry}

\usepackage{graphicx}

\makeatletter
\let\orgdescriptionlabel\descriptionlabel
\renewcommand*{\descriptionlabel}[1]{%
  \let\orglabel\label
  \let\label\@gobble
  \phantomsection
  \edef\@currentlabel{#1}%
  \let\label\orglabel
  \orgdescriptionlabel{#1}%
}
\makeatother

\usepackage[utf8]{inputenc}
\usepackage{latexsym,xspace}
\usepackage{amsmath}
\usepackage{amsfonts}
\usepackage{amsthm}
\usepackage{mathtools}
\newtheorem{lemma}{Lemma}[section]

\theoremstyle{definition}
\newtheorem{definition}{Definition}[section]
\usepackage{amssymb}
\usepackage{bm}
\usepackage[linesnumbered, ruled, noend]{algorithm2e}
\usepackage{lmodern}

\usepackage{ellipsis}

\usepackage{indentfirst}
\usepackage{enumitem}

\usepackage{booktabs}
\usepackage{xcolor}

\usepackage{nameref}
\usepackage{multirow}

\usepackage[numbers,sort]{natbib}
\usepackage{url}

\usepackage{tikz}
\usepackage{lscape}
\usepackage[justification=centering]{caption}
\usepackage{float}
\usepackage{authblk}
\usetikzlibrary{arrows}
\usepackage[pagebackref,pdfdisplaydoctitle,menucolor=orange!40!black,filecolor=magenta!40!black,urlcolor=blue!40!black,linkcolor=red!40!black,citecolor=green!40!black,colorlinks]{hyperref}
\usepackage[textsize=tiny]{todonotes}

\newcommand{\articles}{\ensuremath{{\mathcal{A}}}\xspace}
\newcommand{\article}[1]{\ensuremath{\mathrm{art}({#1})}\xspace}
\newcommand{\location}[1]{\ensuremath{\text{loc}({#1})}\xspace}

\newcommand{\volume}[1]{\ensuremath{\mathrm{vol}({#1})}}
\newcommand{\zones}{\ensuremath{\mathcal{Z}}\xspace}

\newcommand{\zoneOf}[1]{\ensuremath{\text{zone}({#1})}}
\newcommand{\warehouseItems}{\ensuremath{\mathcal{I}}}

\newcommand{\orders}{\ensuremath{\mathcal{O}}\xspace}
\newcommand{\order}{\ensuremath{o}\xspace}
\newcommand{\conveyors}{\ensuremath{\mathcal{C}}\xspace}
\newcommand{\conveyor}[1]{\ensuremath{C_{#1}}}
\newcommand{\distanceOracle}{\ensuremath{d}}

\newcommand{\itemGoal}{\mathit{IG}\xspace}
\newcommand{\containerVolume}{\ensuremath{V}}
\newcommand{\maxOrderPerBatch}{\ensuremath{Q}}

\newcommand{\ordersInBatch}[1]{\ensuremath{\textit{O}_{#1}}}

\newcommand{\picklist}{\ensuremath{p}\xspace}
\newcommand{\cost}[1]{\ensuremath{\text{cost}(#1)}}

\newcommand{\problem}{Joint Order Selection, Allocation, Batching and Picking Problem\xspace}
\newcommand{\problemShort}{JOSABPP\xspace}
\newcommand{\DGA}{DGA\xspace}
\newcommand{\RDGA}{RDGA\xspace}
\newcommand{\selectedItems}{s}
\newcommand{\allSelectedItems}{I}

\makeatletter
\def\moverlay{\mathpalette\mov@rlay}
\def\mov@rlay#1#2{\leavevmode\vtop{%
   \baselineskip\z@skip \lineskiplimit-\maxdimen
   \ialign{\hfil$\m@th#1##$\hfil\cr#2\crcr}}}
\newcommand{\charfusion}[3][\mathord]{
    #1{\ifx#1\mathop\vphantom{#2}\fi
        \mathpalette\mov@rlay{#2\cr#3}
      }
    \ifx#1\mathop\expandafter\displaylimits\fi}
\makeatother

\newcommand{\Perm}[1]{\ensuremath{\mathrm{Perm}({#1})}\xspace}

\title{Joint Order Selection, Allocation, Batching and Picking for Large Scale Warehouses}
\author{{Giorgio Abelli, Maximilian Katzmann, Imran Khan, Olaf Maurer, Julius Pätzold, Pawe\l{} Pszona, Jan-David Salchow}}
\affil{Zalando SE, Berlin, Germany}
\affil[ ]{\it \{giorgio.abelli, max.katzmann, imran.khan, olaf.maurer, julius.paetzold, pawel.pszona, jan-david.salchow\}@zalando.de}
\date{}

\begin{document}
\maketitle

\begin{abstract}
Order picking is the single most cost-intensive activity in picker-to-parts warehouses, and as such
has garnered large interest from the scientific community which led to multiple
problem formulations and a plethora of algorithms published.
Unfortunately, most of them are not applicable at the scale of really large warehouses
like those operated by Zalando, a leading European online fashion retailer.

Based on our experience in operating Zalando's batching system, we propose a novel batching problem
formulation for mixed-shelves, large scale warehouses with zoning. It brings the selection of orders
to be batched into the scope of the problem, making it more realistic while at the same time increasing the optimization potential.

We present two baseline algorithms and compare them on a set of generated instances.
Our results show that first, even a basic greedy algorithm requires significant runtime to solve real-world instances and second, 
including order selection in the studied problem shows large potential for improved solution quality.

\end{abstract}

\section{Introduction}

In manual picker-to-parts warehouses, order picking (the process of collecting items from the warehouse
floor in order to fulfill customer orders) constitutes more than 55\% of the total warehouse cost~\cite{ZHANG2021458}.
On top of that, about 50\% of picker's time is spent traveling between picked items' locations~\cite{tompkins2010facilities}.
Hence, optimizing the process in general, and minimizing its walking component in particular,
have been at the heart of many researchers' interests.

\subsection{Context}

In this section, we describe aspects of warehouse order picking that are relevant
in the setting considered in this paper.
The overall problem entails decisions made on a few levels and, consequently, a few optimization problems
modeling various parts of it have been introduced and studied. We start with an overview of the most prominent ones.
Please refer to the surveys~\cite{de2007design, trendsInPicking, pardo2024} for a more complete picture.

\newpage

\begin{description}
  \item[Picker routing.]
    At the lowest level, items assigned to a picker need to be collected in a single tour.
    As other time components (setup time, searching for items and picking them) are usually considered
    static, travel time is the one aspect which is subject to optimization~\cite{henn2013metaheuristics}.
    As already mentioned, it also accounts for a significant share of the time spent on the picking process,
    so it is only logical to try to optimize it.

    This optimization problem is captured by the Picker Routing Problem (PRP, see~\cite{ratliff1983order, scholz2016new}):
    given a set of item locations a picker should visit, the goal is to produce an ordered {\it picklist}
    which is a complete description of a picker's desired walking route through (a part of) the warehouse.
    A standard approach to optimization of travel time is to consider minimization of the distance
    traveled along the picklist~\cite{THEYS2010755, deKosterPort, hall1993distance}.

  \item[Order batching.]
    To benefit from economies of scale, items belonging to different customer orders are often
    picked together~\cite{de2007design}.
    This is the main idea behind {\it batching} -- grouping customer orders together for joint picking.
    In the general sense, a batch is defined as a logical grouping of a number of orders, together with
    a set of related picklists containing exactly the items assigned to the orders forming the batch.
    The goal of the problem is to create batches covering all orders, with the objective of minimizing
    the sum of lengths of all involved picklists.

    A commonly studied formalization of the problem is known as the Order Batching Problem (OBP),
    (see, e.g., \cite{Waescher2004}). In this publication, a batch is equivalent to a picklist,
    meaning that there is exactly one picklist per batch.
    In other variants of the problem, multiple picklists per batch are permitted~\cite{GV05, yang2020order}.

  \item[Integrated solutions.]
    Given that the objective in the batching problem (OBP) is to minimize total length of picker routes,
    it comes as no surprise that a joint order batching and picking problem was widely studied,
    starting with~\cite{WO05}. Sometimes, other planning aspects such as pick scheduling are included
    in the considered problem as well, e.g., see~\cite{vangils2019, scholz2017}.
    This direction of research, summarized in~\cite{vangils2018} and~\cite{trendsInPicking},
    has produced superior results over approaches treating subproblems of order picking separately.
\end{description}

\noindent We continue by describing how order picking is organized in the large-scale warehouses we study.

\begin{description}
  \item[Mixed-shelves storage.]
    As mentioned in~\cite{boysen2018}, a mixed-shelves storage policy is employed in large scale
    logistics facilities of e-commerce companies the likes of Amazon and Zalando.
    Warehouses operating under this storage policy do not prescribe a specific storage
    location where all items of a particular SKU are to be stored. Instead, items of the same SKU
    can be stored in arbitrary locations throughout the warehouse~\cite{Daniels1998AMF}.

    Such spread happens organically when an unsupervised storage policy is utilized, i.e. when workers freely
    decide where to put the items. Going a step further,
    the {\it explosive storage} policy~\cite{explosivestorage} prescribes a concerted effort
    be taken in order to distribute arriving items of an SKU all over the warehouse.

    Employing a mixed-shelves storage policy introduces another decision to be made in the order picking landscape,
    namely {\it item assignment}: deciding which items and, most importantly, from which locations
    to fulfill an order. Whereas item assignment was considered predetermined and part of the input
    in many classic approaches, recent research~\cite{WEIDINGER2019501, mixedshelvesbatching}
    has shown that incorporating it into the considered problem leads to improved order picking efficiency.

  \item[Zoning.] \label{zoning}
    Even with mixed-shelves storage, distances between items assigned to an order
    can get big in large scale warehouses, making it inefficient to collect
    all of them in a single tour.

    To deal with this problem and improve pick locality, warehouses are subdivided
    into distinct {\it picking zones}. An individual picker is then collecting items only in a zone
    to which they are assigned to. This means that a picklist can only include locations from a single zone.

    When a picklist is finished, the corresponding container is transported to
    a sorting station via a conveyor system. There, items from various zones,
    collected in picklists for the same batch, are consolidated into original
    orders. This is known as {\it pick-and-sort}~\cite{dekoster2012}.
\end{description}

\subsection{Our contribution}

In existing warehouses such as those operated by Zalando, order backlog sizes are orders
of magnitude larger than the few hundred orders which were considered ``large instances'' in the literature.
On top of that, usually there are way more customer orders present than what the warehouse can process
in the nearest future (e.g., in the next hour). This situation brings the following key insight:
instead of batching a set of known orders completely, it is enough to only batch as many
as are expected to be processed during the execution phase we are planning for.

Such an approach has two clear benefits. First, selecting only a subset of orders for batching
opens up optimization potential -- batched orders can be selected in such a way that they complement
each other well when it comes to the creation of efficient picklists.
Second, it is more practical from the operational perspective. While warehouse order picking
is inherently an online problem (new orders and / or items are constantly arriving in the warehouse),
a simplified offline variant is considered instead (only taking orders and items known at a given point in time into account), where the item goal can be used by operations for planning work depending on available worker capacity. This can then be adjusted frequently and in relatively short increments
(as a rule of thumb, between 15 minutes and 2 hours).

We formalize this concept by extending the joint order batching and picker
routing problem (JOBPRP~\cite{JOBPRP}) in the mixed-shelves setting with an
{\it item goal}, i.e., a target number of items to be included in the resulting
batches. Note that we count items instead of orders, since the time required to
process orders is harder to predict due to their varying sizes. Since we are
targeting a {\it zoned warehouse} (see~\ref{zoning}), items belonging to a
batch's picklists need to be separated into orders at a sorting facility.
Additionally, there are physical constraints (like the number of cells in a
manual sort shelf, or
the number of chutes in automated sorters (e.g.,~\cite{linesorters})
which limit the amount of orders that can be processed simultaneously at such a facility.
Therefore, we make the maximum number of orders per batch an additional constraint in our model.

We formalize the above-mentioned problem, which we refer to as \problem (\problemShort, for short), and present a greedy algorithm for solving it.
We conduct an extensive evaluation of its performance using a series of generated instances,
comparing it against a baseline algorithm,
which represents a simplified version of the greedy algorithm.
Our experimental results reveal that the greedy algorithm outperforms the baseline algorithm in terms of the optimization objective.
However, this improvement in performance comes at the cost of increased runtime.
Additionally, we show the impact of batching a subset of orders from a given order pool compared to batching all available orders.
\subsection{Structure of the paper}

The rest of this paper is organized as follows. In Section~\ref{layout-and-terminology}, the targeted
warehouse layout and related nomenclature are introduced.
In Section~\ref{problem-definition}, definition of the considered problem is formalized.
In Section~\ref{greedy-algorithms}, algorithms solving the problem are introduced, and experimental
results of running them are summarized in Section~\ref{experiments}. Section~\ref{conclusion}
summarizes the paper and sketches possible future research directions.

\section{Warehouse layout and terminology} \label{layout-and-terminology}
\begin{figure}[hbt]
    \centering
    \includegraphics{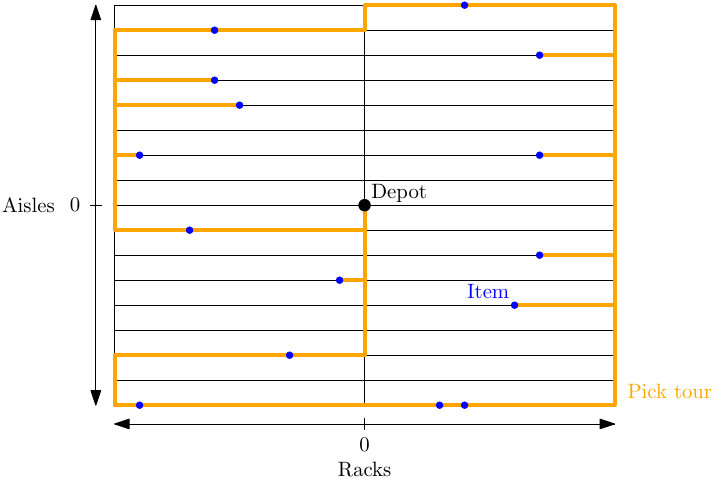}
    \caption{An illustration of a pick tour in a zone, where: the depot is located at the center (black dot),
    items (blue dots) are in aisles and rows, and the orange line depicts the pick tour starting and ending at the depot.}
    \label{fig:zonelayout}
\end{figure}

We start by explaining considered warehouse representation, which we assume is divided into a set
of zones, denoted by \zones.
An individual zone is composed of parallel {\it aisles}, with {\it racks} (shelves used for
storing items) running along them, and a single depot (access point to the conveyor).
Figure~\ref{fig:zonelayout} shows an example pick tour.
It is worth to note that any pick tour must be fully contained in a single zone, as well as begin and finish
at that zone's depot.

\section{Problem definition} \label{problem-definition}

Before giving a formal definition of the \problem, we first give an intuitive overview.

Broadly speaking, the problem input consists of a set of customer orders and a set of available items within one warehouse. Each customer orders a set of articles, potentially containing duplicates.
We assume without loss of generality that there are enough items in the warehouse to satisfy all customer orders.

As mentioned before, the warehouse is divided into zones and each item in the warehouse is located in one of these zones. Together with the zone information, the location of a warehouse item is described by \emph{rack} and \emph{aisle} number. There could be multiple items stored at the same location \nolinebreak-- one can think of a box full of T-shirts, for example.
Each warehouse zone is assumed to have exactly one conveyor station corresponding to a depot, which is the place where warehouse workers pick up and drop off containers.
Hence, a pick tour consists of a warehouse worker picking up a container at the conveyor station, picking all items on the picklist, and returning the full container to the conveyor station. 

To complete the definition of the problem input, it remains to introduce the
previously mentioned constraints. There is a lower bound on the total number of
batched items (the item goal). There is an upper bound on the number of orders
in a batch. The total volume of items in a picklist cannot exceed the volume of
the container used for picking. Note that we assume that a single order never
exceeds this limit.
A formalization of the input is given in Table \ref{tab:definitions}.

\begin{table}[t]
    \small
    \begin{tabular}{ p{0.23\linewidth} p{0.16\linewidth} p{0.5\linewidth}}
        \toprule
        \textbf{Name} & \textbf{Symbol} & \textbf{Description}\\
        \midrule
    Articles & $\articles$ & Each article $a\in\articles$ has a volume $\volume{a}\in\mathbb{R}_{>0}$.\\
    Orders & \orders &  Each order $\order\in\orders$ is a (multi-)set of articles $\{a_1, a_2, \dots\}$.\\
    Zones & $\zones$ & Set of zones\\
    Locations & $\mathcal{L}\subseteq \zones \times \mathbb{Z} \times\mathbb{Z}$ & A location is a triplet of zone, rack and aisle number.\\
    Warehouse Items & $\warehouseItems$ & An item $i\in\warehouseItems$ is associated with an article $\article{i}\in\articles$ and a location $\location{i}\in\mathcal{L}$. For convenience we define $\zoneOf{i} \in \zones$ to be the zone that $i$ is located in.\\
    Conveyor Stations &  \conveyors $\subseteq\mathcal{L}$ &  A set of conveyor stations with exactly one conveyor station per zone $z$ denoted by \conveyor{z}\\
    Walking Distance &$\distanceOracle \colon \mathcal{L} \times \mathcal{L} \to \mathbb{R}_{\geq0}$ & For two locations $x, y$ in different zones, we have $d(x,y) = \infty$.\\
    Item Goal & $\itemGoal \in \mathbb{N}$ & Lower bound on number of items to be batched\\
    Batch Order Limit & $\maxOrderPerBatch \in \mathbb{N}$& Upper bound on number of orders in a batch\\
    Picklist Volume Limit & $\containerVolume \in \mathbb{R}_{>0}$ & Upper bound of total volume of items in a picklist\\
    \bottomrule
\end{tabular}
\caption{Problem Input and Parameters}
\label{tab:definitions}
\end{table}

The objective is to create batches, each batch consisting of a subset of the
orders and a set of picklists, denoting which items to pick to fulfill the
orders, while minimizing the lengths of the corresponding pick tours. Here,
length refers to distance traveled between container pick up and drop off,
which we formalize as the cost of the
corresponding picklist. 

\begin{definition}
A picklist $p=(i_1,\dots,i_{|p|})$ is a sequence of warehouse items $i\in\warehouseItems$ that are all located in the same zone $z \in \zones$. The cost of a picklist consists of walking from the conveyor $C_{z}$ to the first item, visiting all items in order, and walking back to the conveyor. That is,
\[
\cost{p} \coloneqq d(C_{z}, \location{i_1})  + \left(\sum_{j=1}^{|p|-1}d\big(\location{i_j}, \location{i_{j+1}}\big)\right) + d(\location{i_{|p|}}, C_{z}).
\]
\end{definition}

We now turn our attention to the solution space $\mathcal{S}$. A possible
solution $S \in \mathcal{S}$ is a finite set of batches $(O,
P)$, where $O \subseteq \orders$ is called an \emph{order set} and $P$ is a
\emph{picklist set} containing the picklists that are meant to fulfill the
orders in $O$. As described above, a picklist is an ordered subset of items in
the warehouse, i.e., $P \subseteq \{ p \in \Perm{I} \mid I \subseteq
\warehouseItems \}$, where $\Perm{I}$ denotes the set of permutations of $I$.
In the following, we use the shorthand notation $P_S = \bigcup_{(O, P) \in S}
P$ to denote the set of all picklists in a solution $S$.

\begin{definition} For the sake of clarity, we first define the core of the
    \problem, and follow up with additional constraints afterwards.

    The objective (\ref{eq:objective}) is to find a set $S \in \mathcal{S}$
    of batches that minimizes the costs of the picklists, which is subject
    to the following constraints. Each
    order may only appear in one order set or, phrased differently, no two
    order sets may share an order (\ref{eq:empty-order-intersection}). 
    Analogously, no two picklists may share an item (\ref{eq:empty-item-intersection}).
    For a given batch $(O, P)$ the items in
    the picklists in $P$ need to match to the articles in the orders
    in~$O$ (\ref{eq:matching-fulfillment}). Note that $\bigcup$ denotes the additive union of multisets.
    And finally, the item goal needs to be reached (\ref{eq:item-goal}),
    meaning that the total number of items in all picklists is at least as
    large as $\itemGoal$.

\begin{align}
    \min_{S \in \mathcal{S}} \quad & \sum_{p\in P_S} \cost{p} \label{eq:objective}\\
    \text{s.t.} \quad & O \cap O' = \emptyset &\quad \forall (O, P) \neq (O', P') \in S \label{eq:empty-order-intersection}\\
                      & p \cap p' = \emptyset &\quad \forall p \neq p' \in P_S \label{eq:empty-item-intersection}\\
                      & \bigcup_{o \in O} o = \bigcup_{p \in P}\bigcup_{i \in p} \article{i} &\quad \forall (O, P) \in S \label{eq:matching-fulfillment}\\
                      &\sum_{p \in P_S} |p| \ge \itemGoal \label{eq:item-goal},
\end{align}

To complete the problem definition, we add the remaining constraints.
All items in a picklist are in the same zone (\ref{eq:single-zone}).
The total volume of items in a picklist does not
exceed the container volume (\ref{eq:volume-limit}). The number of orders in
a batch does not exceed the prescribed limit (\ref{eq:order-limit}).

\begin{align}
    |\{\zoneOf{i} \mid i\in p\}| &\leq 1 &\quad \forall p\in P_S \label{eq:single-zone}\\
    \sum_{i\in p} \volume{\article{i}} &\leq \containerVolume &\quad \forall p\in P_S \label{eq:volume-limit}\\
    |O| &\leq Q & \quad \forall (O, P) \in S \label{eq:order-limit}.
\end{align}

\end{definition}

\section{Baseline Solution Algorithms} \label{greedy-algorithms}

While the main purpose of this work is to provide a thorough introduction to the \problem, we also want to present baseline solution algorithms.
The goal for these algorithms is not to solve the \problem in the best possible way, but to provide an easy-to-understand entry point for the reader on how a possible solution algorithm could look like.
In this section we propose two algorithms: The \emph{Distance Greedy Algorithm (\DGA)} (Algorithm \ref{alg:one}) 
and a simplified version of it, which we call \emph{Randomized Distance Greedy Algorithm (\RDGA)}.

\DGA, given an instance of the \problem, gradually computes order sets
$O$ and matching picklist sets $P$ until either the item goal is reached or
there are no more orders left in \orders to be processed.
To compute a batch, the algorithm starts with an initially empty set of orders and one by one adds a new order to it together with a set of selected items.
The order is chosen in a way that the corresponding selected items minimize the distance to already selected items (Function \texttt{best\_order}), which is averaged by the number of the selected items.
Orders are added to the set until the stopping criterion is reached, meaning that the number of orders per batch has been reached or there are no orders left.
Once this stopping criterion is met, the selected items of the batch are grouped in their zones and split into picklists.
This is done heuristically: The items are sorted by their aisle number and then clustered together into picklists. The pseudocode is presented in Algorithm \ref{alg:one}.

\begin{algorithm}
	\caption{Distance Greedy Algorithm}\label{alg:one}
	\KwData{as per Table \ref{tab:definitions}}
	$S \gets \emptyset$ \\
	\While {$\itemGoal >0 \; \text{and} \; \orders \neq \emptyset$}{
		$(\ordersInBatch{}, \allSelectedItems) \gets (\emptyset, \emptyset)$ \\
		\While{$|\ordersInBatch{}| < \maxOrderPerBatch \; \text{and} \; \orders \neq \emptyset \; \text{and} \; |\allSelectedItems| < \itemGoal$}{
			$(\order, \selectedItems) \gets \texttt{best\_order}(\orders, \allSelectedItems, \warehouseItems)$  \\
			$(\ordersInBatch{},\; \orders) \gets (\ordersInBatch{}\cup \{\order\},\; \orders \setminus \{\order\})$ \\
			$(\allSelectedItems,\; \warehouseItems) \gets (\allSelectedItems \cup \selectedItems,\;  \warehouseItems \setminus \selectedItems)$\\
		}
		$P \gets \texttt{compute\_picklists}(\allSelectedItems)$ \\
		$\itemGoal \gets \itemGoal - |\allSelectedItems|$\\
		$S \gets S \cup \{(O, P)\}$
	}
	\KwRet $S$\\
	\BlankLine
	\SetKwFunction{FMain}{best\_order}
  \SetKwProg{Fn}{Function}{:}{}
  \Fn{\FMain{\orders, $\allSelectedItems$, \warehouseItems}}{
	$(\order^*, \selectedItems^*, d^*) \gets (null, \emptyset,\infty)$\\
	\For{$\order \in \orders$}
	{
		$d \gets 0$, $\selectedItems \gets \emptyset$ \\
		\For{$a \in \order$} {
			$I_a \gets \{i \in \warehouseItems \mid \article{i} = a \}$\\
			$L \gets \{\location{j} \mid j \in \selectedItems \cup \allSelectedItems  \} \cup \conveyors $\\
			$(i', \ell') \gets \mathrm{argmin}_{i\in \textit{I}_a, \ell \in L} d(\location{i}, \ell) $\\
			$d' \gets d(\location{i'}, \ell')$\\
			$(d,\; \selectedItems) \gets (d + d',\; \selectedItems \cup \{i'\})$ \\
		}
		\If{$\selectedItems^* = \emptyset$ \textbf{or} $\frac{d}{|\selectedItems|} < \frac{d^*}{|\selectedItems^*|}$} {
		    $(\order^*,\;\selectedItems^*,\; d^*) \gets (\order,\;\selectedItems,\;d)$ \\
		}
	}
	\KwRet $(\order^*, \selectedItems^*)$
  }

	\BlankLine
	\SetKwFunction{FPicklist}{compute\_picklists}
	\SetKwProg{Fn}{Function}{:}{}
	\Fn{\FPicklist{$\allSelectedItems$}}{
		$P \gets \emptyset$ \\
		\For{$z \in \zones$}
		{
			$\warehouseItems_z \gets \; \{i \in \allSelectedItems \; | \;\zoneOf{i} =z\}$ \\
			$\picklist \gets \emptyset$ \\
			\For{$i \in \mathrm{sorted}(\warehouseItems_z, \mathrm{key}=(i.\mathrm{aisle}, i.\mathrm{rack}))$}{
				\If{ $\volume{\article{i}} + \sum_{j \in p} \volume{\article{j}} \leq \containerVolume$}{
					$\picklist \gets \picklist \cup \{i\}$ \\
				}
				\Else{
					$ P \gets P \cup \{\picklist\}$\\
					$\picklist \gets \{i\}$\\
				}
			}
			\If{$\picklist \neq \emptyset$}{
				$P \gets P \cup \{\picklist\}$
			}
		}
		
		\KwRet $P$
  }
	
\end{algorithm}
\pagebreak
\begin{lemma}[Complexity of DGA]
    Given a warehouse with $|\warehouseItems| = n$ items, the complexity of Algorithm DGA is $O(n^4)$.
\end{lemma}
\begin{proof}
	To prove the claim, we first show that the running time of Algorithm
	\ref{alg:one} \emph{without} computing picklists is $O(n^4)$ and then
	show that computing \emph{all} picklists takes $O(n \log(n))$ time in total.

	Algorithm \ref{alg:one} makes at most $n$ calls to \texttt{best\_order}
	since every added order has at least one item and the total number of
	requested items $\itemGoal$ is not more than $n$. For the first part it
	thus remains to show that the function \texttt{best\_order} takes
	$O(n^3)$ time. This function iterates through all requested articles of
	all remaining orders, which are in $O(n)$. The minimum operator
	iterates over both $I_a$ and $L$, which has a complexity of $O(n^2)$,
	yielding the desired $O(n^3)$ for \texttt{best\_order}.

	To complete the proof, we show that all calls to
	\texttt{compute\_picklists} take $O(n \log(n))$ time in total. The
	function gets passed an item set $I$, which is split into sets $I_z$
	for the constantly many zones $z \in \zones$. Let $I_1,
	\dots, I_k$ denote all such sets created throughout a run of
	Algorithm \ref{alg:one}, and note that $\sum_{j = 1}^{k} |I_j| \le n$,
	since the number of items in the picklists is bounded by the number of
	items in the warehouse. Computing the sets themselves takes $O(n)$ in
	total and the running time of \texttt{compute\_picklists} is dominated
	by the sorting, which takes time
	\begin{align}
		O\bigg(\sum_{j = 1}^{k} |I_j| \log(|I_j|) \bigg) = O\bigg(\sum_{j = 1}^{k} |I_j| \log(n) \bigg) = O\left(n \log(n) \right). \notag
	\end{align}
	The first equality holds since for all $j$ we have $|I_j| \le n$
	and due to the monotonicity of the logarithm. The second equality
	follows since the number of items in all sets is bounded by the number of
	warehouse items.
\end{proof}

One can see that \texttt{best\_order} dominates the complexity for \DGA by
computing several nested minima. Since this might be very time-consuming we
introduce \RDGA, which works by adjusting \texttt{best\_order} to not iterate
through all orders, but to pick one order at random, while still determining
the best items to fulfill this order. By doing so, we reduce the complexity of
the algorithm from $O(|\warehouseItems|^4)$ to $O(|\warehouseItems|^3 )$. In
the following section, we present our experiments, showing how this modification
significantly improves runtime at the expense of a worsened solution quality.

\section{Experiments} \label{experiments}

In this section, we evaluate the performance of \DGA and \RDGA, for which we compare runtime and solution quality.
We implemented both algorithms (\DGA and \RDGA) in Python~3.10 and made them available at \href{https://github.com/zalandoresearch/batching-benchmarks/}{https://github.com/zalandoresearch/batching-benchmarks/}.
Additionally, we published the instances we used in our experiments and the tool that created them.

\subsection*{Generation of Instances}
The instance generation is based on the problem definition (Section~\ref{problem-definition}) and these assumptions:
\begin{itemize}
    \item Each zone in the warehouse is based on a two-dimensional grid, where one axis (the aisles) is freely walkable,
    and the second axis (the racks) is only walkable along the three cross-aisles.
    The layout is depicted in Figure~\ref{fig:zonelayout}.
    \item Each zone has exactly one depot (modeled as a single node) at $\text{rack}=0$ and $\text{aisle}=0$.
\end{itemize}
These assumptions aim at simplifying the instance format, while at the same time keeping the instances realistic.

We generated 15 instances equally divided into three categories: small, medium
and large, based on the number of orders, items, zones. A zone always has 100
aisles and 100 racks. See Table~\ref{tab:instance-parameters-details} for
details about the parameters for the categories and the individual instances.

\begin{table}
    \centering
    \footnotesize
    \begin{tabular}{llrrclrrclrr}
        \toprule
        Category & \multicolumn{3}{c}{\textbf{Small}} && \multicolumn{3}{c}{\textbf{Medium}} && \multicolumn{3}{c}{\textbf{Large}} \\
        \cmidrule{2-4} \cmidrule{6-8} \cmidrule{10-12}
        Items & \multicolumn{3}{c}{10,000} && \multicolumn{3}{c}{100,000} && \multicolumn{3}{c}{1,000,000} \\
        Orders & \multicolumn{3}{c}{500} && \multicolumn{3}{c}{5,000} && \multicolumn{3}{c}{50,000} \\
        Zones & \multicolumn{3}{c}{10} && \multicolumn{3}{c}{50} && \multicolumn{3}{c}{100} \\
        \cmidrule{2-4} \cmidrule{6-8} \cmidrule{10-12}
        Instances & \textbf{Name} & \textbf{Items} & $\bm{\mathit{IG}}\xspace$ && \textbf{Name} & \textbf{Items} & $\bm{\mathit{IG}}\xspace$ && \textbf{Name} & \textbf{Items} & $\bm{\mathit{IG}}\xspace$ \\
        & \multicolumn{2}{r}{\textbf{\hspace{1.2cm} Demand}} && & \multicolumn{2}{r}{\textbf{\hspace{1.2cm} Demand}} && & \multicolumn{2}{r}{\textbf{\hspace{1.2cm} Demand}} & \\
        \cmidrule{2-4} \cmidrule{6-8} \cmidrule{10-12}
        & small-0 & 1,322 & 264 && medium-0 & 13,115 & 2,623 && large-0 & 131,873 & 26,374 \\
        & small-1 & 1,345 & 269 && medium-1 & 13,223 & 2,644 && large-1 & 131,872 & 26,374 \\
        & small-2 & 1,312 & 262 && medium-2 & 13,135 & 2,627 && large-2 & 131,827 & 26,365 \\
        & small-3 & 1,330 & 266 && medium-3 & 13,236 & 2,647 && large-3 & 131,864 & 26,272 \\
        & small-4 & 1,325 & 265 && medium-4 & 13,135 & 2,625 && large-4 & 132,092 & 26,418 \\
        \bottomrule
    \end{tabular}
    \caption{Parameters for category sizes and instance characteristics. Here, \emph{Items} denotes the number of items in the warehouse and \emph{Items Demand} denotes the sum of all order sizes.}
    \label{tab:instance-parameters-details}
\end{table}

\subsection*{Analysis}
We ran the experiments on a 2020 MacBook Pro with Apple M1 Chip (16GB RAM, 3.2GHz clock rate) on a single thread.
The results are summarized in Tables~\ref{tab:results-greedy} and~\ref{tab:results-random}, respectively.
A comparison of the two tables shows that \DGA yields better results in terms of optimization objective.
On average, \DGA-produced solutions are better than \RDGA solutions by a factor of~2.
This advantage of \DGA, however, goes hand in hand with increased computation times.
More precisely, on large instances \DGA can be up to $1500$ times slower than \RDGA, which takes less than $10$ seconds (as opposed to $9941$ seconds, about 2.5 hours, in the case of \DGA).

\begin{table}[!hbt]
    \centering
    \footnotesize
    \begin{tabular}{lrrrrr}
        \toprule
        \textbf{Instances} & \textbf{Runtime} \textit{(seconds)} & \textbf{Objective Value} & \textbf{Selected Items}& \textbf{Picklists }& \textbf{Batches}\\
        \midrule
        small-0            & 3                                   & 9,122                   & 365                   & 32                         & 3                       \\
        small-1            & 3                                   & 8,796                   & 376                   & 34                         & 3                       \\
        small-2            & 3                                   & 8,372                   & 373                   & 31                         & 3                       \\
        small-3            & 4                                   & 9,546                   & 373                   & 33                         & 3                       \\
        small-4            & 4                                   & 9,226                   & 382                   & 35                         & 3                       \\
        \midrule
        medium-0           & 113                                 & 66,602                  & 2,687                 & 872                        & 23                      \\
        medium-1           & 114                                 & 66,212                  & 2,713                 & 837                        & 23                      \\
        medium-2           & 114                                 & 65,320                  & 2,667                 & 847                        & 23                      \\
        medium-3           & 114                                 & 66,690                  & 2,722                 & 854                        & 23                      \\
        medium-4           & 112                                 & 65,668                  & 2,665                 & 865                        & 23                      \\
        \midrule
        large-0            & 9,972                               & 606,478                 & 26,405                & 10,842                     & 231                     \\
        large-1            & 10,062                              & 607,618                 & 26,434                & 10,830                     & 231                     \\
        large-2            & 9,927                               & 600,110                 & 26,419                & 10,810                     & 231                     \\
        large-3            & 9,955                               & 607,246                 & 26,375                & 10,840                     & 231                     \\
        large-4            & 9,941                               & 610,202                 & 26,440                & 10,907                     & 230                     \\
        \bottomrule
    \end{tabular}
    \caption{Results for \DGA}
    \label{tab:results-greedy}
\end{table}
\begin{table}[!hbt]
    \centering
    \footnotesize
    \begin{tabular}{lrrrrr}
        \toprule
        \textbf{Instances} & \textbf{Runtime} \textit{(seconds)} & \textbf{Objective Value} & \textbf{Selected Items}& \textbf{Picklists}& \textbf{Batches}\\
        \midrule
        small-0            & $<$1                                & 9,814                   & 265                   & 20                         & 2                       \\
        small-1            & $<$1                                & 15,260                  & 398                   & 32                         & 3                       \\
        small-2            & $<$1                                & 10,782                  & 269                   & 20                         & 2                       \\
        small-3            & $<$1                                & 15,334                  & 392                   & 31                         & 3                       \\
        small-4            & $<$1                                & 9,996                   & 273                   & 22                         & 2                       \\
        \midrule
        medium-0           & $<$1                                & 141,752                 & 2,726                 & 951                        & 21                      \\
        medium-1           & $<$1                                & 141,932                 & 2,743                 & 971                        & 21                      \\
        medium-2           & $<$1                                & 137,050                 & 2,634                 & 900                        & 20                      \\
        medium-3           & $<$1                                & 135,246                 & 2,653                 & 905                        & 20                      \\
        medium-4           & $<$1                                & 136,742                & 2,660                 & 897                        & 20                      \\
        \midrule
        large-0            & 7                                   & 1,494,408               & 26,420                & 14,134                     & 201                     \\
        large-1            & 7                                   & 1,489,980               & 26,430                & 14,236                     & 201                     \\
        large-2            & 7                                   & 1,493,316               & 26,401                & 14,255                     & 201                     \\
        large-3            & 7                                   & 1,495,004               & 26,296                & 14,123                     & 200                     \\
        large-4            & 7                                   & 1,497,366               & 26,454                & 14,203                     & 201\\
        \bottomrule
    \end{tabular}
    \caption{Results for \RDGA}
    \label{tab:results-random}
\end{table}

\subsubsection*{Impact of Order Selection}
In a second experiment we investigated the impact of order selection on the overall batch quality. 
With the item goal $\itemGoal$ we defined the \problem to choose a subset and not all of the orders.
It can be seen in Table~\ref{tab:instance-parameters-details} that we are only required to select roughly a fifth of all order items to satisfy the item goal.
Now, we are asking, how does the solution change if we are not given this choice?
To this end, we trimmed down the number of total order articles of the generated instances such that they exactly match the item goal.
More precisely, we reduced the order pool of the input instance by randomly selecting orders until the total count of selected order items exceeded the item goal.
For higher confidence in the results, we repeated this experiment five times for each instance and averaged the objective value of the solution.
Hence, we created five different reduced order pools where the number of order articles equals $\itemGoal$ and reran \DGA on them.

We summarize our findings in Table~\ref{tab:results-comparison-item-goal}. To
be able to compare the quality of solutions of varying sizes, we define the
\emph{picklist cost per item} for a given solution $S$ with picklists $P_S$ as
\begin{align}
    \text{pcpi}(P_S) := \frac{\sum_{p \in P_S} \cost{p}}{\sum_{p \in P_S} |p|}. \tag{picklist cost per item}
\end{align}

    \begin{table}[!hbt]
        \centering
        \footnotesize
        \begin{tabular}{lrrr}
            \toprule
            \textbf{Instances} & \textbf{pcpi Original Instances} & \textbf{pcpi Modified Instances} & \textbf{Difference}\\
            \midrule
            small-0  & 24.99 & 37.76 & +51.1\% \\
            small-1  & 23.39 & 37.62 & +60.8\% \\
            small-2  & 22.45 & 38.36 & +70.9\% \\
            small-3  & 25.59 & 39.13 & +52.9\% \\
            small-4  & 24.15 & 38.06 & +57.6\% \\
            \midrule
            medium-0 & 24.79 & 42.31 & +70.7\% \\
            medium-1 & 24.41 & 42.05 & +72.3\% \\
            medium-2 & 24.49 & 42.28 & +72.6\% \\
            medium-3 & 24.50 & 41.90 & +71.0\% \\
            medium-4 & 24.64 & 41.87 & +69.9\% \\
            \midrule
            large-0  & 22.97 & 39.07 & +70.1\% \\
            large-1  & 22.99 & 39.14 & +70.2\% \\
            large-2  & 22.72 & 39.00 & +71.7\% \\
            large-3  & 23.02 & 39.13 & +70.0\% \\
            large-4  & 23.08 & 39.05 & +69.2\% \\
            \bottomrule
        \end{tabular}
        \caption{Total objective value normalized by number of chosen items}
        \label{tab:results-comparison-item-goal}
    \end{table}

One can observe that the solutions to the new modified instances have a significantly higher picklist cost per item, on average a $66.73\%$ worse solution quality.
Hence, we can conclude that the freedom to choose the orders, i.e. working with a requested number of items $\itemGoal$ leads to significantly better solution quality.
In the literature it is often assumed that all available orders need to be batched (which is the case with the previously defined modified instances), but actually this is not necessarily required.
For example, at Zalando the batching algorithms only need to compute batches such that the warehouse workers have work for the next 30 to 60 minutes, and hence many existing orders are not immediately batched.
This insight from reality together with the computational results from above justify, or even necessitate, the use of an item goal $\itemGoal$.

\section{Conclusion} \label{conclusion}
In this work we presented the \problem that models a central process of Zalando's warehouse operations.
We gave a formal problem definition and justified the need to integrate batching, item allocation and picker routing into one holistic problem.
Furthermore, we explained the need to add an item goal -- the minimum number of requested items -- to the problem definition.

Algorithmically, we presented two baseline approaches to solve the \problem. 
These algorithms intend to serve as a starting point for solving the problem and which can be used for comparisons with more sophisticated algorithms.
The first algorithm, \DGA, chooses the next order for a batch via a heuristic distance-based evaluation of all remaining orders.
The second algorithm, \RDGA, simplifies this approach even further by just selecting a random order as the next order to be added to a batch.

For the computational experiments, we generated a set of instances with parameters chosen to represent real-world applications.
Algorithms, instances and the evaluation procedure were made publicly available.
Based on our experimental results we concluded that the greedy algorithm clearly outperforms the randomized greedy at the cost of higher computation times.
From the runtime explosion of \DGA it is evident how difficult an implementation in reality is since the problem needs to be solved within a limited time frame in order to keep warehouse operations running.
Beyond that, we could see in the second set of experiments that the introduction of an item goal $\itemGoal$ is key to improved solution quality.

One pertinent direction for future research involves presenting more sophisticated algorithms for the \problem, which have to be designed around a careful trade-off between runtime and solution quality.
Another promising direction involves analyzing the optimality gap of the presented solution algorithms.
These can be found by lower bound computations as well as formulating exact approaches.

In general the intention of this publication is to give the interested reader a reasonably low entry barrier into real-world warehouse throughput optimization, in terms of problem complexity and presented baseline algorithms.
Our hope is to make it convenient for other researchers to delve into this problem and we are looking forward to a fruitful exchange of ideas.

{\footnotesize
\bibliographystyle{alpha}
\bibliography{bib}
}

\end{document}